 \documentclass[12pt]{article}
\usepackage{amsfonts}
\usepackage{amssymb}
\usepackage{amsmath}
\usepackage[latin1]{inputenc}

\def\thebibliograph#1#2{\section*{{\normalsize \bf #2}}\list

 {[\arabic{enumi}]}{\settowidth\labelwidth{[#1]}\leftmargin\labelwidth
     \advance\leftmargin\labelsep
     \usecounter{enumi}}
     \def\newblock{\hskip .11em plus .33em minus -.07em}
     \sloppy
     \sfcode`\.=1000\relax}

\newcommand{\re}{{\mathbb{R}}}
\newcommand{\R}{{\mathbb{R}}^n}
\newcommand{\N}{{\mathbb{N}}}

\newtheorem{thm}{Theorem}
\newtheorem{prop}{Proposition}
\newtheorem{DE}{Definition}

\newtheorem{LEM}{Lemma}
\newtheorem{rem}{Remark}

\newcommand{\be}{\begin{equation}}
\newcommand{\ee}{\end{equation}}
\newcommand{\beq}{\begin{eqnarray}}
\newcommand{\eeq}{\end{eqnarray}}
\newcommand{\beqq}{\begin{eqnarray*}}
\newcommand{\eeqq}{\end{eqnarray*}}

{\vskip\baselineskip\noindent\textbf{Proof of Theorem \protect\ref{#1}\quad}}%

{\vskip\baselineskip\noindent\textbf{Proof of Proposition \protect\ref{#1}\quad}}%

{\vskip\baselineskip\noindent\textbf{Proof of Corollary \protect\ref{#1}\quad}}%

{\vskip\baselineskip\noindent\textbf{Proof of Lemma \protect\ref{#1}\quad}}%

\title{ An introduction to composition operators in Sobolev spaces}
  \author{ G\'erard Bourdaud}
 \date{\today}

\begin{document}
\maketitle

\begin{abstract}
We propose a survey of the results on the composition operators in classical Sobolev spaces, obtained between
1975 and 2020. A first version of these notes were the subject of
 a series of lectures, given in Padova University in January 2018.

\end{abstract}



\section{Introduction}

The composition of two maps $f$ and $g$ is defined by $(f\circ g)(x):=f(g(x))$,
if the range of $g$ is included in the definition set of $f$. We denote by $T_f$ the composition operator
$T_f(g):= f\circ g$.

\begin{DE}\label{compdef} Let $E$ be a set of real valued functions, and let $f:\re\rightarrow \re$. We say that
$f$ acts on $E$ by composition (or : superposition) if $T_f(E)\subseteq E$.
\end{DE}

Here are some elementary examples :
\begin{itemize}
\item Let $E$ be a vector space of functions, which means that $g_1+g_2\in E$ and $\lambda g_1\in E$, for all
$g_1,g_2\in E$ and all $\lambda \in\re$. Then every linear function $f:\re\rightarrow \re$ acts on $E$.
\item Let $E$ be an algebra of functions, which means that $E$ is a vector space as above, and that $g_1g_2\in E$  for all
$g_1,g_2\in E$. Then any polynomial $f$ such that $f(0)=0$ acts on $E$.

\end{itemize}

We have a list of natural problems concerning operators $T_f$.

In case $E$ is a vector space of functions, a composition operator $T_f$ is said {\em trivial} if the function $f$ is linear. Then we have the following first question :\\

$\mathcal{Q}_1$: {\em Does exist non trivial composition operators ?}\\

In case $E$ is an algebra of functions, the answer is positive. We will see that it is negative for certain
Sobolev spaces.\\

$\mathcal{Q}_2$: {\em Describe explicitly the set of functions which act on E.}\\

For instance, if $E$ is the set of all continuous functions from $\re$ to $\re$, then a function $f$ acts on $E$ iff 
$f$ is itself continuous.\\

In case $E$ is endowed with a norm, then the following problems make sense :\\

$\mathcal{Q}_3$: {\em Determine the functions $f$ for which $T_f:E\rightarrow E$ is bounded.}\\

$\mathcal{Q}_4$: {\em Determine the functions $f$ for which $T_f:E\rightarrow E$ is continuous.}\\

We propose a wide survey on the answers to the above questions, in case $E$ is the classical Sobolev space $W^m_p(\R)$. Some results are given together with their proofs. Some proofs are simpler than the original ones.

\section{Notation} ${\mathbb N}$ denotes the set of all positive integers, including $0$. ${\mathbb Z}$ denotes the set of all integers. For $x\in\R$, $|x|$ denotes its euclidean norm.\\

If $E,F$ are topological spaces, then $E\hookrightarrow F$ means that $E\subseteq F$, as sets, and the natural mapping $E\rightarrow F$ is continuous. If $B$ is a Lebesgue measurable subset of $\R$, we denotes by $|B|$  its Lebesgue measure. We denote by $\chi_A$ the characteristic function of a set $A$.\\

A {\em multi-index} is $n$-uple $\alpha:=(\alpha_1,\ldots,\alpha_n)\in \N^n$. For such $\alpha$, and for all $h:=(h_1,\ldots,h_n)\in \R$, we set
$|\alpha|:= \alpha_1+\cdots + \alpha_n$ (this differs from the euclidean norm), $\alpha! :=  \alpha_1! \cdots \alpha_n!$,
$h^\alpha := h_1^{\alpha_1}\cdots h_n^{\alpha_n}$. If $f$ is a function defined on an open subset of $\R$, and $\alpha\in \N^n$ as above, we note
$f^{(\alpha)}$ the partial derivative
\[ \frac{\partial^{|\alpha|}f}{\partial x_1^{\alpha_1}\cdots \partial x_n^{\alpha_n}}.\]

If $h\in \R$, the {\em translation operator} is defined by
$(\tau_hf)(x):= f(x+h)-f(x)$ for all function $f$ on $\R$.  The {\em finite difference operator} is defined by $\Delta_hf:=\tau_{-h}f-f$. The $m$-th power of $\Delta_h$ satisfies the following formula :
\be\label{mpower} (\Delta_h^mf)(x)= \sum_{k=0}^m \binom{m}{k} (-1)^{m-k} f(x+kh)\ee
(easy proof by induction).\\

Let $\Omega$ be an open subset of $\R$. We denote by $L_{1,loc}(\Omega)$ the set of (equivalence classes of) locally integrable functions on $\Omega$, endowed with its natural topology (mean convergence on compact subsets of $\Omega$), and by ${\mathcal D}(\Omega)$ the set of all indefinitely differentiable compactly supported functions on 
$\Omega$, endowed with its natural topology, see \cite[1.56]{AF}.\\

Let $Q:=[-1/2,1/2]^n$. We fix some function $\rho\in \mathcal{D}(\R)$ s.t. $\rho(x)=1$ on $Q$ and $\mathrm{supp}\,\rho\subset 2Q$.\\

Let $E$ be a subset of  $L_{1,loc}(\R)$. We say that a function $f\in L_{1,loc}(\R)$ belongs
{\em locally} to $E$ if $\varphi f\in E$ for all
 $\varphi\in {\mathcal D}(\R)$ ; in case $E$ is endowed with a norm, we say that a function $f\in L_{1,loc}(\R)$ belongs
{\em locally uniformly} to $E$ if 
\[ \sup_{a\in \R} \| (\tau_a\varphi)f\|_E<+\infty\,,\]
for all $\varphi\in {\mathcal D}(\R)$.\\

In all the paper, ``ball'' means ``ball with non zero radius'' (we exclude balls reduced to one point).

\section{Composition operators in Lebesgue spaces}

\begin{prop}\label{lp} Let $1\leq p<+\infty$, and let $f:\re\rightarrow \re$ be a Borel function.
Then $f$ acts on $L_p(\R)$ iff there exists $c>0$ s.t.
\be\label{lip0}  |f(t)|\leq c|t|\,,\quad \mathrm{for\,all}\, t\in\re\,.\ee
\end{prop}

{\bf Proof.} 1- If the estimation (\ref{lip0}) holds, it is easily seen that $g\in L_p(\R)$ implies
$f\circ g\in L_p(\R)$. Indeed, the following holds :
\be\label{flp}  \|f\circ g\|_p\leq c\, \|g\|_p,\quad \mathrm{for\,all}\, g\in L_p(\R)\,.\ee

2- Assume that $f$ acts on $L_p(\R)$. Since $L_p(\R)$ does contain nonzero constant functions, it holds $f(0)=0$. Arguing by contradiction, let us assume
that the estimation (\ref{lip0}) does not hold. Then, for some sequence $(a_k)_{k\geq 1}$, we have $|f(a_k)|> k|a_k|$ for all $k\geq 1$. Consider a sequence $(B_k)_{k\geq 1}$ of disjoint measurable sets in $\R$ such that :
\be\label{bk} |a_k|^p |B_k| = k^{-p-1}\,.\ee
Let
\[ g:= \sum_{k\geq 1} a_k\chi_{B_k}\,.\]
By (\ref{bk}), it follows easily that $g\in L_p(\R)$. Since 
\[ f\circ g= \sum_{k\geq 1} f(a_k)\chi_{B_k}\,,\]
 (\ref{bk}) implies again $f\circ g\notin L_p(\R)$, a contradiction.

\begin{rem} {\em  The above proof works as well in case of $L_p(A)$, for any measurable subset $A$ of $\R$ s.t. $|A|=+\infty$. For the generalization of Proposition \ref{lp} to $L_p$ spaces on abstract measure spaces, we refer to \cite[thm.~3.1]{AZ}.}
\end{rem}

In case of linear operators on normed spaces, it is well known that boundedness is equivalent to continuity.
Of course that does not hold for nonlinear ones. In particular, composition operators can be bounded but not continuous.

\begin{prop}\label{cont} Assume $1\leq p\leq +\infty$. Let $(X,\mu)$ be a measure space. Assume that $(X,\mu)$ is non trivial, i.e. there exists a measurable set $A$ in $X$ s.t. $0<\mu(A)<+\infty$. Let $f:\re\rightarrow \re$ be  s.t. $T_f$ takes $L_p(X,\mu)$ to itself.  If  $T_f$ is continuous from $L_p(X,\mu)$ to itself, then $f$ is continuous.
\end{prop}

{\bf Proof.} Assume that $T_f$ is continuous from $L_p(X,\mu)$ to itself. Without loss of generality, assume $f(0)=0$. Let $A$ be  as in the above statement. For all real numbers $u,v$, it holds
\[f\circ  u\chi_A - f\circ v\chi_A= (f(u)-f(v)) \chi_A\,,\]
hence 
\be\label{fuv} \|f\circ  u\chi_A - f\circ v\chi_A\|_p= |f(u)-f(v)| \,\mu(A)^{1/p}\,.\ee
Clearly
\[ \lim_{v\rightarrow u} v\chi_A= u\chi_A\quad \mathrm{in}\,\, L_p\,.\]
By continuity of $T_f$, and by (\ref{fuv}), we obtain the continuity of $f$.\\

By Propositions \ref{lp} and \ref{cont}, it follows that, in case of $L_p(\R)$, there exist bounded composition operators which are not continuous.  Proposition \ref{cont}  admits a reciprocal :

\begin{prop}\label{cont+} Let $f:\re \rightarrow \re$ be a continuous function s.t., for some constant
$c>0$, it holds $|f(t)|\leq c\,|t|$, for all $t\in \re$. Let $(X,\mu)$ be a measure space and let $1\leq p<+\infty$. Then $T_f$ is continuous from $L_p(X,\mu)$ to itself.
\end{prop}

{\bf Proof.} It suffices to prove the following : for all sequence $(g_j)$ converging to $g$ in $L_p(X,\mu)$, there exists a subsequence  $(g_{j_k})$ s.t. $(f\circ g_{j_k})$ converges to $f\circ g$ in
$L_p(X,\mu)$.
By a classical measure theoretic result (see, for instance, \cite{R2}), there exists a subsequence  $(g_{j_k})$ and a function $
h\in L_p(X,\mu)$ s.t. 
\[g_{j_k} \rightarrow g\quad a.e. \,,\quad \left|g_{j_k} \right| \leq h\,.\]
By continuity of $f$, it holds $f\circ g_{j_k} \rightarrow  f\circ g$ a.e..
By assumption on $f$, it holds
\[ \left |f\circ g_{j_k} -f\circ g\right|\leq 2c\,h\,.\]
By Lebesgue dominated convergence Theorem, we conclude that $\left \|f\circ g_{j_k} -f\circ g\right\|_p$ tends to $0$.

\begin{rem}{\em  If $f: \re \rightarrow \re$ is bounded and continuous, $T_f$ is easily seen to be continuous from $L_\infty(X,\mu)$ to itself. The details are left to the reader.
}\end{rem}

\section{Automatic boundedness}

\begin{DE}\label{bound}
Let $E$ be a normed space. A mapping $T:E\rightarrow E$ is said {\em bounded} if, for all bounded set $A$ of $E$, the set $T(A)$ is bounded.
\end{DE}

For instance, according to the estimation (\ref{flp}), any composition operator, which takes $L_p(\R)$ to itself, is bounded on $L_p(\R)$. More generally,
for all ``reasonable'' function space, a weak form of boundedness is satisfied by composition operators.
Thus we have a kind of automatic boundedness for a large class of function spaces.

\begin{prop}\label{weakb} Let $E,F$ be vector  subspaces of $L_{1,loc}(\Omega)$. Assume that :
\begin{itemize}
\item $E$ and $F$ are endowed with complete norms s.t. the embeddings of $E$ and $F$ into
$ L_{1,loc}(\Omega)$ are continuous.
\item ${\mathcal D}(\Omega)$ is embedded  into $E$.
\item For all $\varphi \in{\mathcal D}(\Omega)$ and $g\in F$, it holds $\varphi g\in F$.
\end{itemize}
For all $f:\re\rightarrow \re$ s.t. $f(0)=0$ and $T_f(E)\subseteq F$, there exist a closed ball $B\subset \Omega$ and two numbers
$c_1,c_2>0$ such that, for all $g\in E$,
\be\label{weakb1}\|g\|_E\leq c_1 \quad \mathrm{and} \quad \mathrm{supp}\, g\subseteq B\qquad \Rightarrow \qquad \|f\circ g\|_F\leq c_2\,.\ee
\end{prop}

{\bf Proof.} By contradiction, assume that, for all $B,c_1,c_2$ there exists $g\in E$ s.t.
\be\label{contra} \|g\|_E\leq c_1 \,,\quad \mathrm{supp}\, g\subseteq B\,,\quad \|f\circ g\|_F>c_2\,.\ee
Consider a sequence $(B_j)_{j\geq 1}$ of disjoint closed balls in $\Omega$. Take functions $\varphi_j\in  {\mathcal D}(\Omega)$ s.t
$\varphi_j(x)=1$ on $\frac{1}{2} B_j$ (the ball of same center and half radius than $B_j$) and $\varphi_j(x)=0$ out of $B_j$. It is easily seen (Closed Graph Theorem, see \cite[chap.~II, §6, thm.~1]{Y} or \cite[thm.~2.15]{R1}) that, for $\varphi \in  {\mathcal D}(\Omega)$, the linear multiplication operator $g\mapsto \varphi g$ is bounded on $F$. Thus we can consider
\[M_j:= \sup\{ \| \varphi_j g\|_F\,:\, \|g\|_F\leq 1\}\,.\]
According to (\ref{contra}), there exist functions $g_j$ s.t.
\[ \|g_j\|_E\leq 2^{-j} \,, \quad \mathrm{supp}\, g_j\subseteq \frac{1}{2}B_j\,,\quad \|f\circ g_j\|_F>jM_j\,.\]
Let $g:= \sum_j g_j$. Clearly $g\in E$ and, by embedding $E\hookrightarrow L_{1,loc}(\Omega)$, it holds
\[ g(x)= \sum_{j\geq 0} g_j(x) \quad \mathrm{a.e.}\,.\]
By considering supports, it holds $\varphi_j(f\circ g)= f\circ g_j$, hence
\[ jM_j<  \|\varphi_j(f\circ g)\|_F\leq M_j \|f\circ g\|_F\]
for all $j\geq 1$, a contradiction.

\begin{rem} {\em If $\Omega = \R$, and if $E$ is translation and dilation invariant, the conclusion of Proposition \ref{weakb} can be improved : indeed {\em for all} ball or cube $B$, there exists $c_1,c_2>0$ s.t. (\ref{weakb1}) holds for all $g\in E$.
}\end{rem}

As an example of use of Proposition \ref{weakb}, we give the following variant of Proposition \ref{lp} :

\begin{prop}\label{lpb} Let $1\leq p<+\infty$, let $\Omega$ be an open subset of $\R$ s.t. $|\Omega|<+\infty$, and let $f:\re\rightarrow \re$ be a Borel function.
Then $f$ acts on $L_p(\Omega)$ iff there exists $\alpha,\beta>0$ s.t.
\be\label{lip0b}  |f(t)|\leq \alpha |t|+\beta\,,\quad \mathrm{for\,all}\,\, t\in\re\,.\ee
\end{prop}

{\bf Proof.} Since the sufficiency of condition (\ref{lip0b}) is clear, we deal only with necessity.
Assume that $f$ acts on $L_p(\Omega)$. Without loss of generality, we can assume that $f(0)=0$. By Proposition \ref{weakb}, there exist a cube $Q'\subset \Omega$ and two numbers
$c_1,c_2>0$ such that, for all $g\in L_p(\Omega)$,
\be\label{weakb2}\|g\|_p\leq c_1 \quad \mathrm{and} \quad \mathrm{supp}\, g\subseteq Q'\qquad \Rightarrow \qquad \|f\circ g\|_p\leq c_2\,.\ee
It  holds $Q'= b+ 2rQ$ for some $b\in \Omega$ and $r>0$. For any $a\in \re$, and $0<\varepsilon\leq 1$, let
\[ g_{a,\varepsilon}(x):= a \rho\left(  \frac{x-b}{r \varepsilon}  \right)\,.\]
Then the support of $ g_{a,\varepsilon}$ is included into $Q'$. Concerning the norm of  $g_{a,\varepsilon}$, we have to discuss according to $a$.

In case of large $a$, more precisely if $|a|\geq R:=r^{-n/p}c_1 \|\rho\|_p^{-1}$, we choose $\varepsilon =
\varepsilon (a)$ s.t.
\be\label{adjust} |a|\,r^{n/p} \|\rho\|_p\,\varepsilon^{n/p}=c_1\,.\ee

If $|a|<R$, we take $\varepsilon=1$. 
In both cases, we obtain $\|g_{a,\varepsilon}\|_p\leq c_1$, hence 
$\|f\circ g_{a,\varepsilon}\|_p\leq c_2$.
Since \[\rho\left(  \frac{x-b}{r\varepsilon}\right)=1\]
on the cube $b+\varepsilon rQ$, this implies

\[ \int_{b+\varepsilon rQ} |f(a)|^p\,\mathrm{d}x\leq c_2^p\,,\]
hence $|f(a)|^p \varepsilon^n\leq c_3$, for some constant $c_3$.

If $|a|\geq R$, by using (\ref{adjust}), we obtain $|f(a)|\leq c_4|a|$, for some constant $c_4$.
If $|a|<R$, we obtain $|f(a)|\leq c_3^{1/p}$. This ends up the proof.

\begin{rem} {\em The above proof can be viewed as a prototype of a number of results on composition operators, as we will see further.}
\end{rem}

\section{Definition and main properties of Sobolev spaces}

\begin{DE}\label{defwd} Let $f\in L_{1,loc}(\R)$, and $\alpha\in \N^n$. We say that $f$ admits a weak derivative of order $\alpha$ if there exists $g\in L_{1,loc}(\R)$ s.t.
\[ \int_{\R} g(x)\varphi (x)\,{\mathrm d}x=(-1)^{|\alpha|}\int_{\R} f(x)\varphi^{(\alpha)} (x)\,{\mathrm d}x\]
for all $\varphi\in 
{\mathcal D}(\R)$.
\end{DE}

If such $g$ exists, it is easily seen to be unique, up to equality almost everywhere ; then we denote it by $f^{(\alpha)}$ and we call it the {\em weak derivative} of $f$ at order $\alpha$.
\begin{DE}\label{defsob} Let $m\in \N$ and $1\leq p \leq +\infty$. The Sobolev space $W^m_p(\R)$ is the set of functions 
$f\in L_{1,loc}(\R)$ s.t., for all $|\alpha|\leq m$, $f^{(\alpha)}$ exists in the weak sense, and $f^{(\alpha)}\in L_p(\R)$.
\end{DE}

$W^m_p(\R)$ is a vector subspace of  $L_p(\R)$. It will be endowed with the following norm :
\be\label{norm} \|f\|_{W^m_p(\R)}:= \sum_{|\alpha|\leq m} \|f^{(\alpha)}\|_p\,.\ee

We give here the useful properties of Sobolev spaces. 

First of all, $W^m_p(\R)$ is a function space which satisfies the assumptions of Proposition \ref{weakb},
see \cite[thm.~3.3]{AF}.

The behavior of the norm (\ref{norm}) w.r.t. dilations is described in the following assertion, with easy proof :

\begin{prop}\label{dilation} It holds
\[ \| f(\lambda(.))\|_{W^m_p(\R)}\leq  \lambda^{m-(n/p)} \| f\|_{W^m_p(\R)}\,,\]
for all $\lambda\geq 1$.
\end{prop}

Then we have the so-called {\em Sobolev embedding theorems}, see \cite[thm.~4.12]{AF} :

\begin{prop}\label{sobemb}
If
\[ m_1-m_2\geq \frac{n}{p_1}- \frac{n}{p_2}> 0\,.\]
then $W^{m_1}_{p_1}(\R)\hookrightarrow W^{m_2}_{p_2}(\R)$.
\end{prop}

In particular
$W^{m}_{p}(\R)\hookrightarrow L_\infty(\R)$ if $m>n/p$. Indeed we have a more precise statement, where
$C_b(\R)$ denotes the space of bounded continuous functions on $\R$:

\begin{prop}\label{sobemb1} If $ m>n/p$, or $p=1$ and $m=n$, it holds
$ W^{m}_{p}(\R)\hookrightarrow C_b(\R)$.
\end{prop}

On the contrary,
$ W^{m}_{p}$ is not embedded in  $L_\infty(\R)$ in case $m<n/p$, or $m=n/p$ and $p>1$.

\begin{rem} {\em The members of $W^{m}_{p}(\R)$ are  equivalence classes of functions w.r.t. to a.e. equality. Thus the precise meaning of Proposition \ref{sobemb1} is the following :  all $f\in W^{m}_{p}(\R)$ admits a (necessarily unique) bounded continuous representative s.t. $\|f\|_\infty \leq c \|f\|_{W^m_p}$, for some constant $c>0$ depending only on $m,p,n$.}

\end{rem}

\begin{prop}\label{alg}
If $m>n/p$, or $m=n$ and $p=1$, the Sobolev space $W^{m}_{p}(\R)$ is a subalgebra of $C_b(\R)$.
\end{prop}

See  \cite[thm.~4.39]{AF} for the proof.

\section{Necessity of Lipschitz continuity}

In case $m\geq 1$, any function which acts on $W^m_p$ by composition is necessarily Lipschitz continuous, at least locally. This is a major difference w.r.t. to the case of $L_p$. To prove this property, some preliminary results are needed :

\begin{LEM}\label{bump} Assume that $W^m_p(\R)$ is not embedded into $L_\infty(\R)$. There exists a sequence $(\theta_j)_{j\geq 1}$ in
 $\mathcal{D}(\R)$ s.t.
 \[ \theta_j(x)=1\quad \mathrm{on}\quad 2^{-j}Q\,,\quad \mathrm{supp}\,\theta_j\subseteq Q\,,\quad \lim_{j\rightarrow +\infty} \|\theta_j\|_{W^m_p(\R)}=0\,.\]
 \end{LEM}
 
 {\bf Proof.} 
 In case $m<n/p$, we take $\theta_j(x):=\rho(2^jx)$. By Proposition \ref{dilation}, it holds
 \[ \|\theta_j\|_{W^m_p(\R)}\leq 2^{j(m-(n/p))} \|\rho \|_{W^m_p(\R)}\,.\]
 Thus the sequence $(\theta_j)$ has the wished properties.
 
 Now assume that $m=n/p$ and $1<p<+\infty$. Let
 \[\theta_j(x):=\frac{1}{j}\sum_{k=1}^{j}\rho(2^kx)\,.\]
If  $|\alpha|=m$, then the function $x\mapsto \rho^{(\alpha)}(2^kx)$ is supported by the set  $S_k:= 2^{-k+1} Q\setminus 2^{-k}Q$.
Thus, for all $1\leq k\leq j$ and $ x\in S_k$, it holds
  \[\left|\theta_j^{(\alpha)}(x)\right|=\frac{1}{j}2^{mk}\left|\rho^{(\alpha)}(2^kx)\right|\leq cj^{-1}2^{mk}\,.\]
Hence  
\[ \|\theta_j^{(\alpha)}\|_p^p = \sum_{k=1}^{j}\int_{S_k} \left|\theta_j^{(\alpha)}(x)\right|^p\,\mathrm{d}x \leq cj^{-p} \sum_{k=1}^{j} 2^{kmp}2^{-nk}= cj^{1-p}\,.\]
Thus  the sequence $(\|\theta_j^{(\alpha)}\|_p)$ tends to $0$ for all $|\alpha|=m$. The same holds, with easy proof, for $|\alpha|<m$.
That ends up the proof of Lemma \ref{bump}.

\begin{LEM}\label{findif1} Define the sequence of functions $(B_m)_{m\geq 1}$  in $L_1(\re)$ by $B_1:={\bf 1}_{[0,1]}$ and $B_{m+1}:= B_m\ast B_1$ for all $m$.
Then
\be\label{findif2} \Delta^m_hf(x)= \int_{-\infty}^{+\infty} B_m(t)\left( \sum_{|\alpha|=m} \frac{m!}{\alpha!} f^{(\alpha)}(x+th)\,h^\alpha\right)\mathrm{d}t\,,\ee
for almost all $x\in \R$, all $h\in \R$, all $m\geq 1$ and all $f\in W^m_p(\R)$.
\end{LEM}

{\bf Proof.} We deal with a continuously $m$ times differentiable function $f$. An approximation procedure would give the general case.

{\em Step 1 : case $n=1$.} In such case, formula (\ref{findif2}) reduces to
\be\label{findif3} \Delta^m_hf(x)= \int_{-\infty}^{+\infty} B_m(t)  f^{(m)}(x+th)\,h^m\mathrm{d}t\,.\ee
We prove it by induction. The case $m=1$ is well known. Assuming that (\ref{findif3}) holds, we obtain
\[\Delta^{m+1}_hf(x)= h^m\int_{-\infty}^{+\infty} B_m(t)  \Delta_hf^{(m)}(x+th)\,\mathrm{d}t =\]
\[=h^{m+1}\int_{-\infty}^{+\infty} B_m(t)\left( \int_{-\infty}^{+\infty}B_1(s) f^{(m+1)}(x+(t+s)h)\,\mathrm{d}s\right)
\mathrm{d}t
\,.\]
By Fubini, a change of variable, and the definition of $B_m$, we obtain formula (\ref{findif3}) at rank $m+1$.

{\em Step 2 : general case.} We fix $x,h$ in $\R$, and set $g(t):= f(x + t(h/|h)))$ for all $t\in \re$. It holds $\Delta_h^mf(x)= \Delta_{|h|}^mg(0)$.
Then applying Step 1 to the function $g$, we obtain (\ref{findif2}). The details are left to the reader.

\begin{LEM}\label{findif} For all $m\geq 1$, $1\leq p\leq \infty$, there exists $c>0$ s.t.
\[ \left( \int_{\R} \left|\Delta^m_hf(x)\right|^p\,\mathrm{d}x\right)^{1/p}\leq c |h|^m\|f\|_{W^m_p(\R)}\]
for all $h\in \R$ and all $f\in W^m_p(\R)$.

\end{LEM}

{\bf Proof.} By definition of $B_m$, it holds $B_m\geq 0$ and $\int_{-\infty}^{+\infty} B_m(t)\,\mathrm{d}t =1$.
Applying (\ref{findif2}), we obtain
\[ \| \Delta^m_hf \|_p \leq |h|^m  \sum_{|\alpha|=m} \frac{m!}{\alpha!} \|f^{(\alpha)}\|_p\,,\]
hence the desired result.

\begin{thm}\label{Lip} Assume that $m\geq 1$ and that $W^m_p(\R)$ is not embedded into $L_\infty(\R)$. Then any function $f:\re\rightarrow \re$, s.t. $T_f$ takes $W^m_p(\R)$ to itself, is Lipschitz continuous on $\re$.
\end{thm}

{\bf Proof.}
 In all the proof, $\|-\|$ will denote the norm in $W^m_p(\R)$.\\

{\em Step 1 : construction of the comb-shaped function.} This construction was first introduced by S.~Igari \cite{Ig}.
Let $\displaystyle{
A_N :={\mathbb Z}^n\cap[-N,N]^n}$,  for all positive integer $N$. We fix a real number $s$ s.t.
\be\label{sss} 0<s<\frac{1}{2m+1}\,.\ee

Let $b,b'$ be distinct real numbers. Then we consider integers $N,j\geq 1$, and a real number $r>0$, whose values will be fixed below  depending on $b,b'$. Our test function will be defined by
 \begin{equation}\label{thetest}
g(x):=  \sum_{\mu\in A_N} \rho\left( \frac{1}{s}\left( \frac{x}{r}-\mu\right)\right)(b'-b) + \theta_j (x)\,b\,,\end{equation}
 see Lemma \ref{bump} for the construction of $\theta_j$.
The first condition on parameters will be
\begin{equation}\label{rNnu} 3rN \leq 2^{- j}\,.\end{equation}
 By inequality $s <1/2$ and by condition (\ref{rNnu}), we deduce that the cubes 
$r(2s Q+\mu)$, $\mu\in {\mathbb Z}^n$, are pairwise disjoint, and that
$r(2s Q+\mu)\subset r(Q+\mu)\subset 2^{-j}Q$, for all $\mu\in A_N$. Hence
 \begin{equation}
\label{value1} 
g(x) = b'\,, \quad \mathrm{if}\,\,x \in r(s Q+\mu)\,\,\mathrm{for\, some}\,\,
 \mu\in A_N\,,
\end{equation}
\begin{equation}
\label{value2}
g(x) = b\,, \quad\mathrm{if}\,\,x \in 2^{-j}Q\setminus 
\bigcup_{ \mu\in A_N} r(2s
Q+\mu)\,.
\end{equation}
By (\ref{rNnu}), it holds $r\leq 1$. By Proposition \ref{dilation}, we deduce \begin{equation}
\label{atomic}
  \Big\| \sum_{\mu\in A_N}\rho \Big( \frac{1}{s}
 \Big( \frac{{\bf .}}{r}-\mu \Big) \Big)  \Big\|\leq c_1r^{(n/p)-m} N^{n/p}\,, 
\end{equation}
for some constant $c_1$.\\

{\em Step 2 : adjustment of parameters.} Now we assume that $f$ acts on $W^m_p(\R)$ by composition.
By Proposition \ref{weakb}, we can find
constants $\delta_1,\delta_2$  such that $\|f\circ u\| \leq \delta_2$ for all function $u$ s.t. $\|u\|\leq \delta_1$, and $u$ is supported by $Q$. In order to apply this property to $u=g$, we need the following inequalities :
\begin{equation}\label{bnu} |b|\,\|\theta_j\| \leq \frac{\delta_1}{2}\,,\end{equation}
\begin{equation}\label{Nrbb'}
 \frac{\delta_1}{3c_1|b-b'|} \leq r^{(n/p)-m} N^{n/p}\leq
\frac{\delta_1}{2c_1|b-b'|}\,.\end{equation}
 Now we discuss the choice of $j, N, r$ with respect to $b,b'$, s.t. conditions (\ref{rNnu}), (\ref{bnu}) and (\ref{Nrbb'}) hold. First, we choose $j = j(b)\geq 1$ such that (\ref{bnu}) holds. This is possible by Lemma \ref{bump}.
 Then we have to discuss according to $m<n/p$ or $m=n/p$.

In case $m<n/p$, we define
\[ r:= \left( \frac{ \delta_1}{ 2c_1|b-b'|}N^{-n/p}\right)^{ \frac{p}{n-mp}}\,,\] which ensures condition (\ref{Nrbb'}) ;
since \[ rN = \left( \frac{ \delta_1}{ 2c_1|b-b'|}\right)^{ \frac{p}{n-mp}} N^{ \frac{mp}{pm-n}}\,,\]
the condition (\ref{rNnu}) holds  for all sufficiently large $N$, depending on $|b-b'|$.

In case $m=n/p$, 
 we take $N$ s.t. (\ref{Nrbb'}) holds.
Such a choice is possible if $|b-b'|\leq c_2$, where $c_2>0$ depends only on $p,n,\delta_1$. Then we put
$r :=  2^{-j}/3N$.\\

{\em Step 3 : end of the proof.} By combining inequalities (\ref{atomic}), (\ref{bnu}) and (\ref{Nrbb'}), we deduce $ \|g\|\leq \delta_1$. Using Lemma \ref{findif}, we deduce   
\[ \| \Delta^m_h(f \circ g)\|_p \leq  \delta_3|h|^m\,,\]
for all $h\in \R$,
for some constant $\delta_3$ depending only on $\delta_2,m,n,p$. Let $Q^+:=]0,1/2]^n$ and $e_1:=(1,0, \ldots,0)\in \R$. By condition (\ref{sss})
 we have
\[ 
x +\ell rs e_1 \in  r(Q+\mu)\subset2^{-j}Q \,\quad (\ell =0,\ldots , m)  \,,\]
\[
x +\ell rs e_1\notin 
\bigcup_{ \mu'\in A_N} r(2s
Q+\mu')\,,\quad (\ell=1,\ldots , m) \,,
\]
for  all $x \in r(s Q^++\mu)$; for such $x$, equalities (\ref{value1}) and (\ref{value2}), 
and formula  (\ref{mpower}), imply :
\[ 
\left|\Delta^{m}_{rs e_1}(f\circ g)(x)\right| = |f(b')- f(b)|\,.
\]
Hence
\[ \delta_3\geq
 c_{3}r^{-m} \left( \sum_{\mu\in A_N}
\int_{r(s Q^++\mu)} 
\left|\Delta^{m}_{rs e_1}(f\circ g)(x)\right|^p{\rm d}x\right)^{1/p}
\]
\[ \geq c_4 |f(b')- f(b)| \,N^{n/p}r^{(n/p)-m} \,.
\]
By (\ref{Nrbb'}) we obtain the existence of a constant $\delta_4$
s.t.
 $|f(b')-f(b)|\leq \delta_4|b-b'|$ for all $b,b'\in \re$ satisfying
 $|b'-b|\leq c_2$. Thus $f$ is uniformly Lipschitz continuous.

 \begin{thm}\label{Liploc} Assume that $m\geq 1$. Then any function $f:\re\rightarrow \re$, s.t. $T_f$ takes $W^m_p(\R)$ to itself, is locally Lipschitz continuous on $\re$.
\end{thm}

{\bf Proof.} Let $f:\re\rightarrow \re$ be a function which acts on $W^m_p(\R)$.
Let $a\in \re$. We introduce a localized version of $T_f$ with the help of the following statement :

\begin{LEM}\label{local} Under the above assumptions, there exists a nonlinear operator $U_a$ which takes $W^m_p(\R)$ to itself, s.t, for all $g\in {W^m_p(\R)}$ :
\[ U_ag(x)= f(a+g(x))-f(a)\,,\quad \mathrm{for\,all}\quad x\in Q\,,\]
\[\|g\|_{W^m_p(\R)}\leq \delta_1\quad \mathrm{and} \quad \mathrm{supp}\, g\subseteq Q\qquad \Rightarrow \qquad \|U_ag\|_{W^m_p(\R)}\leq \delta_2\,.\]

\end{LEM}

The proof is similar to that of Proposition \ref{weakb}, see \cite[Lemma 1]{BL} for details.\\

Returning to the proof of Theorem \ref{Liploc},
we argue in the same way as in the proof of Theorem \ref{Lip}, just replacing $T_f$ by $U_a$. We define $g$ by (\ref{thetest}), with $\theta_j(x)$ replaced by $\rho(2x)$, $s=1/4$ and $r=1/6N$.
The inequality (\ref{bnu}) becomes $|b|\leq \delta_3$, for some constant $\delta_3$ depending only on $\delta_1$.
The  double inequality (\ref{Nrbb'}) reduces to
\begin{equation}\label{goodN} 
 \frac{\delta_4}{|b-b'|} \leq  N^m\leq \frac{\delta_5}{|b-b'|}\,,\end{equation}
for some constants $\delta_4,\delta_5$ depending on $\delta_1$ and $c_1$. If $|b-b'|\leq \delta_4$, we can choose $N$ satisfying (\ref{goodN}).
We obtain a constant $\delta_6$ s.t.
\[|f(a+b)- f(a+b')|\leq \delta_6|b-b'|\,,\] for all $b,b'$ satisfying $|b|\leq \delta_3$ and $|b-b'|\leq \delta_4$. Thus $f$ is Lipschitz continuous in a neighborhood of $a$.\\

An easy modification of the above proof gives us the following statement :

\begin{prop}\label{holloc} Let us assume $m\geq 3$, and 
 define $p_1$ by :
\begin{equation}\label{otherp} 2-\frac{n}{p_1} := m- \frac{n}{p}\,.\end{equation}
Then every function $f:\re\rightarrow \re$, such that $T_f$ takes  $W^m_p(\R)$ to  $W^2_{p_1}(\R)$,
is locally H\"older continuous with exponent $2/m$.\end{prop}

\section{A case of degeneracy : Dahlberg Theorem}

As announced in Introduction, Sobolev spaces provide simple examples of spaces for which the answer to question $\mathcal{Q}_1$ is negative.

\begin{thm}\label{dahl}  Assume that $m$ is an integer satisfying
\begin{equation}\label{triv} 1+ \frac{1}{p} < m < \frac{n}{p}\,.\end{equation}
Then, for all function $f:\re\rightarrow \re$ which acts on $W^m_p(\R)$ by composition, there exists $c\in \re$ s.t. $f(t)=ct$ for all $t\in \re$.
\end{thm}

This theorem was first proved by B.~Dahlberg \cite{DA}, assuming $f$ of class $C^\infty$. 
Indeed, a slightly stronger property holds true :

\begin{prop}\label{strongdahl} Under condition (\ref{triv}), let us define $p_1$ by condition (\ref{otherp}).
Then, for all function $f:\re\rightarrow \re$ such that $T_f$ takes  $W^m_p(\R)$ to  $W^2_{p_1}(\R)$, there exists $c\in \re$ s.t. $f(t)=ct$ for all $t\in \re$.
\end{prop}

By Proposition \ref{sobemb}, the embedding $W^{m}_{p}(\R)\hookrightarrow W^{2}_{p_1}(\R)$
holds true. Thus Theorem \ref{dahl} follows by Proposition \ref{strongdahl}.\\

{\bf Proof of Prop. \ref{strongdahl}.}  

{\em Step 1.} We assume first that $f$ of class $C^2$. 
Since $W^m_p(\R)$ does not contain nonzero constant  functions, we have $f(0)=0$.  By Proposition \ref{weakb}, there exist two numbers
$c_1,c_2>0$ such that, for all $g\in W^m_p(\R)$,
\be\label{dahl1} \|g\|_{W^m_p(\R)}\leq c_1 \quad \mathrm{and} \quad \mathrm{supp}\, g\subseteq 2Q\qquad \Rightarrow \qquad \|f\circ g\|_{W^2_{p_1}(\R)}\leq c_2\,.\ee
Define the function $u\in {\mathcal D}(\R)$ by
\be\label{functionu} u(x):=x_1\rho(x)\,,\ee
where $x_1$ denotes the first coordinate of $x\in \R$.
Let $a>0$, and $0<\varepsilon \leq 1$
(a number to be determined w.r.t. $a$). Let us define $g_a\in {\mathcal D}(\R)$ by
\[ g_a(x):= a u\left(\frac{x}{\varepsilon}\right)\,.\]
Then $ \mathrm{supp}\, g_a\subset 2Q$, and $\|g_a\|_{W^m_p(\R)}\leq c_1$ if
\be\label{epsilon}  a \,\varepsilon^{(n/p)-m}\|u\|_{W^m_p(\R)} = c_1\,.\ee
Due to assumption $m<n/p$, the above equality determines $\varepsilon$ as a function of $a$, if $a$ is sufficiently large. Hence it holds
$\|f\circ g_a\|_{W^2_{p_1}(\R)}\leq c_2$ for all large $a$'s. Since \[(f\circ g_a)(x)= f\left(\frac{a}{\varepsilon}x_1\right)\,,\quad x\in \varepsilon Q\,,\]
we deduce
\[ \left(\frac{a}{\varepsilon}\right)^{2p_1}\int_{\varepsilon Q} \left|   f''\left(\frac{a}{\varepsilon}x_1\right)  \right|^{p_1}\,{\mathrm d}x\leq c_2^{p_1}\,.\]
By using (\ref{epsilon}) and a change of variable, we obtain a constant $c_3>0$ s.t.
\be\label{last} a^{p_1-1} \int_{-a/2}^{+a/2}|   f''(t)  |^{p_1}\,{\mathrm d}t\leq c_3\,,\ee
for all large $a$'s. By assumption $m>1+(1/p)$, it holds $p_1>1$. If we take $a$ to $+\infty$, we deduce
\[\int_{-\infty}^{+\infty}\left|   f''(t)   \right|^{p_1}\,{\mathrm d}t=0\,\]
Hence $f''(t)=0$ for almost every $t\in \re$. Since $f''$ is continuous, we conclude that $f(t)= ct$, for some constant $c$.\\

{\em Step 2.} We turn now to the general case. By Theorem \ref{Liploc} and Proposition \ref{holloc}, we know that $f$ is continuous. 
Let $\omega\in \mathcal{D}(\re)$, with support in $[-1,+1]$, even, s.t. $\int\omega (t)\,\mathrm{d}t =1$. Let us set $\omega_j(t):=j\omega (jt)$ for all positive integer $j$. The convolution $\omega_j\ast f$ is defined, and it is a smooth function. Let us define
\[ f_j(t):= (\omega_j\ast f)(t) - (\omega_j\ast f)(0)\,.\]
For all function $g$ with support in $Q$, it holds
\[ (f_j\circ g)(x) = \rho(x)\, \int_\re \left( f( (g(x) +t)\rho (x)) - f(t\rho (x))\right) \omega_j(t)\,\mathrm{d}t \]
for all $x\in \R$. In other words :
\be\label{dahl2} \mathrm{supp}\,g\subseteq Q \quad
\Rightarrow \quad
f_j\circ g = \rho\, \int_\re \left(f \circ ( (g+t)\rho ) - f\circ (t\rho) \right) \omega_j(t)\,\mathrm{d}t \,.\ee

Let $M:= \sup\{ \|\rho\,h\|_{W^m_p(\R)}\,:\, \|h\|_{W^m_p(\R)}\leq 1\}$.
Let $j_0$ be the first integer s.t.
\[ j_0\geq  2 c_1^{-1}\|\rho\|_{W^m_p(\R)} \,.\]

Let $g$ be s.t.  $\mathrm{supp}\,g\subseteq Q$ and
 \[\|g\|_{W^m_p(\R)} \leq \frac{c_1}{2M}\,.\]
 Then, for all $j\geq j_0$, and all $|t|\leq 1/j$, it holds
 \[  \|(g+t)\rho\|_{W^m_p(\R)} \leq c_1\,.\]
By (\ref{dahl1}), we obtain
\[   \|f_j\circ g\|_{W^2_{p_1}(\R)} \leq 2Mc_2\,\]
for all $j\geq j_0$. All together, we have obtained constants $c_3,c_4>0$ s.t.
\be\label{dahl3} \|g\|_{W^m_p(\R)}\leq c_3 \quad \mathrm{and} \quad \mathrm{supp}\, g\subseteq Q\qquad \Rightarrow \qquad \|f_j\circ g\|_{W^2_{p_1}(\R)}\leq c_4\,,\ee
for all $j\geq j_0$. Arguing exactly as in the first part of the proof,  we conclude that, for some constants $a_j$, $j\geq j_0$, we have $f_j(t)=a_jt$ for all $t\in \re$. Thus we obtain:
\[(\omega_j\ast f)(t) = (\omega_j\ast f)(0) + a_jt\] for all $t\in \re$.
Since $f$ is continuous, we know that $\lim_{j\rightarrow +\infty}(\omega_j\ast f)(t) = f(t)$
for all $t\in \re$. Taking $t=1$, we obtain $\lim_{j\rightarrow +\infty}a_j= f(1)$.
We conclude that $f(t)=f(1)t$ for all $t\in \re$.

\section{Composition operators on $W^1_p$}First of all, we recall a classical result :

\begin{thm}\label{rade} For all  $f:\re\rightarrow \re$, the three following properties are equivalent :
\begin{itemize}
\item[(1)] $f$ is Lipschitz continuous,
\item[(2)] $f$ admit a weak derivative in $L_\infty(\re)$,
\item[(3)] There exists $g\in L_\infty(\re)$ and a constant $c\in \re$ s.t.
\[ \forall x\in \re\quad  f(x) = \int_0^x g(t)\,\mathrm{ d}t + c\,.\]
\end{itemize}
\end{thm}

{\bf Proof.} The implication $(3)\Rightarrow (1)$ is immediate. The equivalence $(2)\Leftrightarrow (3)$ is easy to prove.
Concerning $(1)\Rightarrow (3)$, we refer to \cite[thm.~7.18]{F} (Alternatively, we can observe that any Lipschitz continuous function is absolutely continuous, then apply \cite[thm.~8.17]{R2}).\\

\begin{thm}\label{w1p} Let $f:\re\rightarrow \re$, s.t. $f(0)=0$. Then $f$ acts on $W^1_p(\R)$ iff
\begin{itemize}
\item $f$ is Lipschitz continuous, in case $W^1_p(\R)\not\subset L_\infty(\re)$,
\item $f$ is locally Lipschitz continuous, in case $W^1_p(\R)\subset L_\infty(\re)$.
\end{itemize} 
\end{thm}

{\bf Proof.} This theorem is due to Marcus and Mizel \cite{MM}. Roughly speaking, sufficiency result relies upon the formula $\partial_j(f\circ g) = (f'\circ g)\partial_jg$. In case $W^1_p(\R)\subset L_\infty(\R)$, we just need that $f'$ belongs to $L_\infty$ on the range of $g$. The necessity of Lipschitz conditions follows by Theorems \ref{Lip} and \ref{Liploc}.
 \section{Full description of acting functions in higher order Sobolev spaces}

Let us give first a sufficient condition for composition :

\begin{thm}\label{w2p} Assume that 
$m\geq \max(2,n/p)$, or $m=2$, $p=1$. If a function $f:\re\rightarrow \re$ satisfies $f(0)=0$ and
$f'\in W^{m-1}_p(\re)$,
then $f$ acts on $W^m_p(\R)$.
\end{thm}

{\bf Proof.} A preliminary remark : under assumption of Theorem \ref{w2p}, it holds $W^{m-1}_p(\re)\hookrightarrow L_\infty(\re)$.
That follows by Proposition \ref{sobemb1}. 

Here we limit ourselves to the case $m=2$. The method that we use is typical of the general case.
Also we assume that $f$ is of class $C^m$, with bounded derivatives up to order $m$, and that $g$ is smooth, with derivatives tending to $0$ at infinity ;
see \cite{B,Bou_09a} and \cite[5.2.4, thm.~2]{RS} for the approximation procedure to cover the general case.

Let $g\in W^2_p(\R)$.
We have to prove that the second order derivatives of $f\circ g$ belongs to $L_p(\R)$.
It holds
\be\label{w2p1} \partial_j\partial_k(f\circ g) = (f''\circ g) (\partial_jg)(\partial_kg) +  (f'\circ g) \partial_j\partial_kg\,.\ee
The second term belongs to $L_p$, because $f'\in L_\infty$. Thus we can concentrate on the first one. By applying Cauchy-Schwarz inequality, we obtain

\be\label{w2p2} \|(f''\circ g) \partial_jg\,\partial_kg\|_p \leq U_j^{1/2p} U_k^{1/2p}\,,\ee where
\[U_j:= \int_{\R} |(f''\circ g)(x)|^p | \partial_jg(x)|^{2p}\,\mathrm{d}x\,.\]

Let us introduce
\[ h(x):= \int_x^{+\infty} |f''(t)|^p\,\mathrm{d}t\,.\]
Then
\[U_j = - 
\int_{\R} (h'\circ g)(x) \partial_jg(x)\partial_jg(x) | \partial_jg(x)|^{2p-2}\,\mathrm{d}x 
= -\int_{\R} \partial_j(h\circ g)(x)\,\partial_jg(x) | \partial_jg(x)|^{2p-2}\,\mathrm{d}x \,.\]
An IP gives
\[U_j= (2p-1)  \int_{\R} (h\circ g)(x)\,\partial_j^2g(x) | \partial_jg(x)|^{2p-2}\,\mathrm{d}x \,.\]
Hence
\be\label{w2p3} U_j\leq (2p-1) \|f''\|^p_p \int_{\R} |\partial_j^2g(x)|\, | \partial_jg(x)|^{2p-2}\,\mathrm{d}x \,.\ee
In case $p=1$, the above inequality becomes
$U_j \leq \|f''\|_1\, \|\partial_j^2g\|_1$. That ends up the proof of Theorem  in case $m=2$, $p=1$.

In case $p>1$, we use H\"older inequality to derive
\[ U_j\leq  (2p-1)\|f''\|^p_p\,\|\partial_j^2g\|_p \left(\int_{\R} | \partial_jg(x)|^{2p}\,\mathrm{d}x \right)^{1-(1/p)}\,.\]
By applying Proposition \ref{sobemb} and condition $2\geq n/p$, it holds
$W^2_p(\R)\hookrightarrow W^1_{2p}(\R)$. That ends up the proof of Theorem \ref{w2p}.

\begin{rem}{\em   The above proof shows also that the composition operator is bounded, under assumptions of Theorem \ref{w2p}.
More precisely, there exist a constant $c=c(p,n)>0$ such that
\be\label{boundw2p}
\|f\circ g\|_{W^2_p(\R)} \leq c  \|f''\|_p\,
\left (\|g\|_{W^2_p(\R)}+\|g\|_{W^2_p(\R)}^{2-(1/p)}\right)\,.\ee
 }\end{rem}

We turn now to the complete description of composition operators. Due to Theorems \ref{dahl} and \ref{w1p}, we will consider only the case $m\geq 2$, together with the three following subcases :
  \begin{itemize}
  \item  $m >n/p$, or $m=n$ and $p=1$.
  \item $m =n/p$ and $p>1$.\item
  $m=2$, $p=1$ and $n\geq 3$ .
 \end{itemize}

   \begin{thm}\label{m>n/p} Let $m\geq 2$, $1\leq p<+\infty$. If $m>n/p$, or if $m=n$ and $p=1$, then a function  $f:\re\rightarrow \re$ acts on $W^m_p(\R)$ iff $f(0)=0$ and $f$ belongs locally to $W^m_p(\re)$.
    \end{thm}
    
    {\bf Proof.} 1- Assume that $f$ belongs locally to $W^m_p(\re)$, and that $g\in W^m_p(\R)$. By Proposition \ref{sobemb1}, $g$ is bounded.
    Let $\varphi\in \mathcal{D}(\re)$ s.t. $\varphi(t)=1$ on the range of $g$. Then $f\circ g=
    (\varphi\,f)\circ g$. Since $\varphi\,f\in W^m_p(\re)$, we can apply Theorem \ref{w2p}, and conclude that $f\circ g\in W^m_p(\R)$.
    
    2- Assume that $T_f$ takes $W^m_p(\R)$ to itself. By considering $f\circ g$, where $g\in \mathcal{D}(\R)$ satisfies $g(x)=x_1$ on an arbitrary ball of $\R$, we conclude that $f$, together with all its derivatives up to order $m$, belong to $L^p$ on each bounded interval of $\re$.
    
    \begin{thm}\label{m=n/p} Let $m =n/p\geq 2$ and $p>1$. Then a function  $f:\re\rightarrow \re$ acts on $W^m_p(\R)$ iff $f(0)=0$ and $f'$ belongs locally uniformly to $W^{m-1}_p(\re)$.
    \end{thm}
    
    {\bf Proof.} The sufficiency of the condition on $f$ follows by a modification of the proof of Theorem \ref{w2p}, see \cite{B,Bou_09a} or \cite[5.2.4, thm.~2]{RS}.
    
    To prove the necessity, we use the same ideas as in the proof of Theorem \ref{Lip}. Let $f:\re\rightarrow \re$ be a function which acts on $W^m_p(\R)$.
We introduce constants $\delta_1,\delta_2$ as in the proof of Theorem \ref{Lip}.
Let $b$ be a real number. Let $j=j(b)\geq 1$ s.t.  (\ref{bnu}) holds. Let us consider the function :
\[g_{b}(x):= \lambda u(2^jx) + \theta_j (x)\, b
\,,
\]
 where $u$ is the function introduced in (\ref{functionu}), and $\lambda$ is a constant, to be fixed below. 
 By assumption $m=n/p$, it holds $\|u(2^j(.))\|\leq \|u\|$. Thus, the choice of $\lambda:= \delta_1/2\|u\|$
 implies $\|g_b\|\leq \delta_1$.
 Hence we have
    \begin{equation}\label{normfg1} \|f\circ g_{b}\|\leq \delta_2\,.\end{equation}
    On the cube $2^{-j}Q$, it holds $
(f\circ g_{b})(x) = f( \lambda 2^{j}x_1 +b)$, hence
\[ 
\partial_1^m(f\circ g_{b})(x) = \lambda^m2^{jm} f^{(m)}( \lambda 2^{j}x_1 +b)\,.
\]
Then using (\ref{normfg1}),  a change of variable,  and condition $m=n/p$, we find a constant $\delta_3>0$ such that
\[  \int_{b-(\lambda/2)}^{b+(\lambda/2)} |f^{(m)}(y)|^p\,{\rm d}y \leq \delta_3 \,,\]
whatever be $b\in \re$. Thus we have proved that $f^{(m)}$ belongs to $L_p(\re)$ locally uniformly.
Since we know yet that $f'\in L_\infty$, it follows easily that $f'$ belongs to $W^{m-1}_p(\re)$ locally uniformly.

     \begin{thm}\label{w21} If $n\geq 3$, then a function  $f:\re\rightarrow \re$ acts on $W^2_1(\R)$ iff $f(0)=0$ and $f''\in L_1(\re)$.
    \end{thm}
    
    {\bf Proof.} Sufficiency of $f''\in L_1(\re)$ follows by Theorem \ref{w2p}. To prove necessity, we proceed so as in the proof of Theorem
    \ref{dahl}. Then the estimation (\ref{last}) becomes
    \[  \int_{-a/2}^{+a/2}|   f''(t)  |\,{\mathrm d}t\leq c_3\,,\]
    for all large $a$. By taking $a\rightarrow +\infty$, we obtain $f''\in L_1(\re)$.\\

\section{Continuity of composition on Sobolev spaces}

The more precise versions of Theorems \ref{w2p},\ref{m>n/p},\ref{m=n/p} show that all composition operators which take
 $W^m_p(\R)$ to itself are bounded. They are also continuous, according to the following :

\begin{thm}\label{contsob}  Let $m$ be an integer $\geq 1$, $1\leq p<\infty$, and let $f:\re\rightarrow \re$.
If $f$ acts by composition on $W^m_p(\R)$, then the composition operator $T_f$ is continuous from
 $W^m_p(\R)$ to itself.
 \end{thm}
 
 This theorem was proved step by step between 1976 and 2019 :
 \begin{itemize}
 \item For $m=1$ and $p=2$, by Ancona \cite{AA},
  \item For $m=1$ and any $p$, by Marcus and Mizel \cite{MM2},
   \item For $m>n/p$ and $1<p<\infty$, by Lanza de Cristoforis and the author \cite{BL},
   \item In the general case by Moussai and the author \cite{BM19}, who proved also this ``automatic'' continuity on the so-called Adams-Frazier spaces 
   $W^m_p\cap \dot{W}^1_{mp}(\R)$, where $\dot{W}$ denotes the homogeneous Sobolev space, and on the spaces $\dot{W}^m_p\cap \dot{W}^1_{mp}(\R)$, conveniently realized.
 \end{itemize}

G\'erard Bourdaud
 
Universit\'e Paris Cit\'e, I.M.J. - P.R.G (UMR 7586)

B\^atiment Sophie Germain
       
Case 7012  
    
75205 Paris Cedex 13     \\                          
bourdaud@math.univ-paris-diderot.fr

\end{document}